\documentclass[a4paper,11pt,french]{smfart}

\usepackage{amssymb}
\usepackage{url}
\usepackage[T1]{fontenc}
\usepackage[cal=boondoxo]{mathalfa}

\makeatletter
\def\MakeTitle{\par \@topnum\z@
  \begingroup
  \let\@makefnmark\relax  \let\@thefnmark\relax
  \@MakeTitle
  \@endMakeTitlehook
  \endgroup
  \c@footnote\z@
  \let\MakeTitle\relax \let\@MakeTitle\relax }
\def\@endMakeTitlehook{}
\def\smf@journalhead{}
\def\@MakeTitle{\cleardoublepage\thispagestyle{copyright}
 \begingroup
 \uppercasenonmath\shorttitle
 \ifx\@empty\shortauthors \let\shortauthors\shorttitle
  \else \uppercasenonmath\shortauthors \andify\shortauthors \fi
  \toks@\@xp{\shortauthors}\@temptokena\@xp{\shorttitle}%
  \edef\@tempa{\@nx\markboth{\the\toks@}{\the\@temptokena}}\@tempa
 \topskip\z@skip
 \vtop to 55 mm{%
 \parindent=0pt
 \hrule
 \medskip
 {\abstractfont\smf@journalhead\par}\vfil
 \begin{center}
 \def\baselinestretch{1.2}\large\vfil
   {\bfseries\smf@boldmath\MakeUppercase\titre\par}
 \vfil
   \ifx\@empty\smfbyname\else
     {\smfbyfont\smfbyname\ifsmf@byauthor\par\vfil\else\ \fi}%
   \fi
   {\edef\smfandname{{\noexpand\normalfont \smfandname}}
    \andify\authors\auteur\par}
 \vfil \vrule height .4pt width .3\textwidth \vfil
 \end{center}}
 \@MakeTitlehook
 \par\bigskip
 \ifx\@empty\@dedicatory\else\@setdedicatory\medskip\fi
 \@setabstract\par\smallskip\@setaltabstract\par
 \bigskip\bigskip
 \endgroup}
\newif\ifsmf@byauthor\smf@byauthortrue
\def\smfbyfont{\normalfont\itshape}
\def\@MakeTitlehook{%
 \ifx\@empty\@subjclass\else\@footnotetext{\@setsubjclass}\fi
 \ifx\@empty\@keywords\else\@footnotetext{\@setkeywords}\fi
 \ifx\@empty\thankses\else\def\par{\let\par\@par}\@footnotetext{\@setthanks}\fi
}
\makeatother
\usepackage{fancybox}
\usepackage{makecell}

\usepackage[leftbars]{changebar}
\usepackage{todonotes}

\usepackage{palatino}
\usepackage{rotating}
\usepackage{graphicx}

\usepackage{enumerate}
\usepackage{enumitem}

\usepackage{color}
\definecolor{shadecolor}{gray}{0.90}

\def\proofname{D\'emonstration}

\input xypic
\xyoption{all}
\xyoption{arc}

\makeindex

\def\indexname{Index des notations}
\newtheorem{theo}{Th\'eor\`eme}[section]
\def\thetheo{\thesection.\arabic{theo}}

\def\theprop{\thesection.\arabic{prop}}
\newtheorem{lem}[theo]{Lemme}
\def\thelem{\thesection.\arabic{lem}}
\newtheorem{coro}[theo]{Corollaire}
\def\thecoro{\thesection.\arabic{coro}}

\renewcommand\theequation{\thetheo}
\def\equat{\refstepcounter{theo}\begin{equation}}
\def\endequat{\end{equation}}

\renewcommand\thetable{\thetheo}

\renewcommand\thesection{\arabic{section}}
\renewcommand\thesubsection{\thesection.\Alph{subsection}}

\def\contentsname{Table des mati\`eres}
\def\listfigurename{Liste des figures}
\def\listtablename{Liste des tables}
\def\appendixname{Appendice}

\def\refname{R\'ef\'erences}

    \def\CM{{\mathbb{C}}}

    \def\FM{{\mathbb{F}}}

\def\Cb{{\mathbf C}}

\def\Gb{{\mathbf G}}

\def\Lb{{\mathbf L}}

\def\Ob{{\mathbf O}}    
  \def\pb{{\mathbf p}}  
    
    \def\RC{{\mathcal{R}}}
\def\Sb{{\mathbf S}}

\def\Drm{{\mathrm{D}}}    
\def\Erm{{\mathrm{E}}}

\def\Zrm{{\mathrm{Z}}}

  \def\gti{{\tilde{g}}}

  \def\wti{{\tilde{w}}}

\def\b{\beta}
\def\g{\gamma}

\def\e{\varepsilon}

\def\l{\lambda}
\def\L{\Lambda}

\def\o{\omega}
\def\O{\Omega}

\def\z{\zeta}

\def\mub{{\boldsymbol{\mu}}}

\DeclareMathOperator{\Id}{{\mathrm{Id}}}

\DeclareMathOperator{\Ker}{{\mathrm{Ker}}}

\def\to{\rightarrow}
\def\longto{\longrightarrow}

\def\fonction#1#2#3#4#5{\begin{array}{rccc}
{#1} : & {#2} & \longto & {#3}  \\
& {#4} & \longmapsto & {#5} 
\end{array}}

\def\lexp#1#2{\kern\scriptspace\vphantom{#2}^{#1}\kern-\scriptspace#2}

\def\ge{\hspace{0.1em}\mathop{\geqslant}\nolimits\hspace{0.1em}}

\mathchardef\inferieur="321E
\mathchardef\superieur="321F

\def\eqna{\begin{eqnarray*}}
\def\endeqna{\end{eqnarray*}}

\def\itemth#1{\item[${\mathrm{(#1)}}$]}

\catcode`\@=11
\long\def\@car#1#2\@nil{#1}
\long\def\@first#1#2{#1}
\long\def\@second#1#2{#2}
\long\def\ifempty#1{\expandafter\ifx\@car#1@\@nil @\@empty
  \expandafter\@first\else\expandafter\@second\fi}
\catcode`\@=12

\newcommand{\REF}{{\mathrm{Ref}}}

\theoremstyle{remark}
\newtheorem{rema}[theo]{Remarque}

\theoremstyle{plain}

\def\BIL{LR}
\def\GAUCHE{L}
\def\CAR{CAR}
\def\FAM{FAM}

\def\xyinj{\ar@{^{(}->}}
\def\xysur{\ar@{->>}}

\bigskip

\makeatletter
\def\hlinewd#1{%
\noalign{\ifnum0=`}\fi\hrule \@height #1 %
\futurelet\reserved@a\@xhline}
\makeatother

\newlength\epaisLigne

\usepackage{multirow}

\makeatletter
\def\hlinewd#1{%
\noalign{\ifnum0=`}\fi\hrule \@height #1 %
\futurelet\reserved@a\@xhline}
\makeatother

\usepackage{multirow}

\usepackage{lscape}

\begin{document}

%\baselineskip=16pt
%\large\baselineskip=20pt
%\Large\baselineskip=24pt

\title{Une construction du groupe ${\boldsymbol{G_{32}}}$ de Shephard-Todd \\
via le groupe de Weyl de type ${\boldsymbol{E_6}}$}

\author{{\sc C\'edric Bonnaf\'e}}
\address{IMAG, Universit\'e de Montpellier, CNRS, Montpellier, France} 

\makeatletter
\email{cedric.bonnafe@umontpellier.fr}
\makeatother

\date{18 octobre 2023}

%\thanks{The author is partly supported by the ANR:
%Projects No ANR-16-CE40-0010-01 (GeRepMod) and ANR-18-CE40-0024-02 (CATORE).}

% 

%\tableofcontents

% \noindent{\bf Acknowledgements.} We thank warmly D. Juteau 
% for providing us the main ideas leading to Example~\ref{ex:b2-g2}.
% 
% \bigskip

\begin{abstract}
Il est connu que le quotient du groupe d\'eriv\'e du groupe de r\'eflexions complexe
$G_{32}$ de Shephard-Todd
(qui est de rang $4$) par son centre est isomorphe au groupe d\'eriv\'e du groupe de
Weyl de type $E_6$. Nous montrons que cet isomorphisme peut se r\'ealiser via la
puissance ext\'erieure deuxi\`eme et en profitons pour proposer une construction
alternative du groupe $G_{32}$.
\end{abstract}

\maketitle

\pagestyle{myheadings}

\markboth{\sc C. Bonnaf\'e}{\sc Une construction de $G_{32}$ via $E_6$}

\bigskip

Notons $G_{32}$ le groupe de r\'eflexions complexe construit par Shephard-Todd~\cite{ST} et
notons $E_6$ un groupe de Weyl de type $\Erm_6$ (qui est not\'e $G_{35}$ dans la classification
de Shephard-Todd). Notons $\Sb\pb_4(\FM_{\! 3})$ (resp. $\Sb\Ob_5(\FM_{\! 3})$)
le groupe symplectique (resp. orthogonal) de dimension $4$ (resp. $5$) sur le corps fini \`a trois
\'el\'ements $\FM_{\! 3}$. Notons $\O_5(\FM_{\! 3})$ l'image de $\Sb\pb_4(\FM_{\! 3})$
dans $\Sb\Ob_5(\FM_{\! 3})$ \`a travers le morphisme naturel
$\Sb\pb_4(\FM_{\! 3}) \to \Sb\Ob_5(\FM_{\! 3})$~: c'est le sous-groupe distingu\'e d'indice $2$
de $\Sb\Ob_5(\FM_{\! 3})$. Pour finir, si $G$ est un groupe, notons respectivement $\Drm(G)$ et $\Zrm(G)$
son groupe d\'eriv\'e et son centre et, si $d$ est un entier naturel non nul,
notons $\mub_d$ le groupe des racines $d$-i\`emes de l'unit\'e dans $\CM$.

Il est montr\'e dans~\cite[th\'eo.~8.43~et~8.54]{lehrer taylor}
que $G_{32} \simeq \mub_3 \times \Sb\pb_4(\FM_{\! 3})$
et que $E_6 \simeq \Sb\Ob_5(\FM_{\! 3})$. En particulier,
$$\Drm(G_{32})/\mub_2 \simeq \O_5(\FM_{\! 3}) \simeq \Drm(E_6).\leqno{(*)}$$
L'objet de cette note est de pr\'esenter une explication \'el\'ementaire directe
de l'isomorphisme $\Drm(G_{32})/\mub_2 \simeq \Drm(E_6)$, qui permet en fait de construire
le groupe de r\'eflexions complexe $G_{32}$ \`a partir du groupe de r\'eflexions rationnel $E_6$.
Cette construction utilise le morphisme classique $\Sb\Lb_4(\CM) \longto \Sb\Ob_6(\CM)$, et suit
les m\^emes lignes que notre pr\'ec\'edent article~\cite{g31} (dans lequel nous construisions
le groupe de r\'eflexions complexe $G_{31}$ en partant du groupe de Weyl de type $B_6$).

Cette note ne pr\'etend pas d\'emontrer un r\'esultat profond~: nous la voyons plut\^ot comme
un joli exemple d'application de la th\'eorie classique des groupes de r\'eflexions
(invariants, th\'eorie de Springer,...).

\bigskip

\noindent{\bf Remarque.} Il est montr\'e dans~\cite[th\'eo.~8.53]{lehrer taylor} que
$G_{33} \simeq \mub_2 \times \O_5(\FM_{\! 3})$. En particulier, cela donne indirectement
un isomorphisme $\Drm(G_{33}) \simeq \Drm(E_6)$. Mais nous ne connaissons pas de construction
de cet isomorphisme qui soit dans le m\^eme esprit que ci-dessus.

\bigskip

\section{Le morphisme $\Sb\Lb_4(\CM) \to \Sb\Ob_6(\CM)$}

\medskip

Nous rappelons ici la construction de ce morphisme, en reprenant certaines notations de~\cite{g31}.
Fixons un espace vectoriel complexe $V$ de dimension $4$ et notons
$$\fonction{\L}{\Gb\Lb(V)}{\Gb\Lb(\wedge^2 V)}{g}{\wedge^2 g}$$
le morphisme naturel de groupes alg\'ebriques. Notons que $\wedge^2 V$ est de dimension $6$
et que
\equat\label{eq:noyau}
\Ker \L \simeq \mub_2 = \{\pm \Id_V\}.
\endequat
Les listes de valeurs propres d'\'el\'ements de $\Gb\Lb(V)$ ou $\Gb\Lb(\wedge^2 V)$
seront toujours donn\'ees avec multiplicit\'es.
Si $g \in \Gb\Lb(V)$ admet pour valeurs propres
$a$, $b$, $c$ et $d$, alors
\equat\label{eq:vp}
\text{\it $\L(g)$ admet $ab$, $ac$, $ad$, $bc$, $bd$ et $cd$ comme valeurs propres.}
\endequat
En particulier,
\equat\label{eq:det}
\det \L(g) = (\det g)^3.
\endequat

Fixons maintenant un g\'en\'erateur $\e$
de l'espace vectoriel unidimensionel $\wedge^4 V$. Le choix de ce g\'en\'erateur permet
d'identifier $\CM$ et $\wedge^4 V$ et de d\'efinir la forme bilin\'eaire
$$\fonction{\b_\wedge}{\wedge^2 V \times \wedge^2 V}{\CM}{(x,y)}{x \wedge y.}$$
Cette forme bilin\'eaire est sym\'etrique et non d\'eg\'en\'er\'ee. Par d\'efinition du
d\'eterminant,
\equat\label{eq:beta}
\b_\wedge(\L(g)(x),\L(g)(y))=(\det g) \b_\wedge(x,y)
\endequat
pour tous $g \in \Gb\Lb(V)$ et $x$, $y \in \wedge^2 V$. Pour des raisons de dimension et de connexit\'e,
l'image de $\Gb\Lb(V)$ \`a travers $\L$ est la composante neutre $\Cb\Ob(\wedge^2 V)^\circ$
du groupe conforme orthogonal
$\Cb\Ob(\wedge^2 V)=\Cb\Ob(\wedge^2 V, \b_\wedge)$ et $\L$ induit un isomorphisme de groupes alg\'ebriques
\equat\label{eq:iso}
\Sb\Lb(V)/\mub_2 \simeq \Sb\Ob(\wedge^2 V).
\endequat

Le lemme \'el\'ementaire suivant nous sera utile
(dans cet article, nous noterons $j$ une racine primitive cubique de l'unit\'e).

\bigskip

\begin{lem}\label{lem:j}
Soit $g \in \Sb\Ob(\wedge^2 V)$ admettant pour valeurs propres $j$, $j$, $j$, $j^2$, $j^2$ et $j^2$.
Alors il existe un unique \'el\'ement $\gti \in \Sb\Lb(V)$, admettant
$1$, $j$, $j$ et $j$ comme valeurs propres, et tel que
$$\L^{-1}(\{g,g^{-1}\}) = \{\pm \gti,\pm \gti^{-1}\}.$$
\end{lem}

\bigskip

\begin{proof}
Soit $h \in \Sb\Lb(V)$ tel que $\L(h)=g$. Alors $\L^{-1}(\{g,g^{-1}\}) = \{\pm h,\pm h^{-1}\}$.
Notons $a$, $b$, $c$ et $d$ les valeurs propres de $h$.
D'apr\`es~\eqref{eq:vp}, quitte \`a r\'eordonner les valeurs propres
de $h$, on peut supposer que $ab=ac=j$. En particulier, $b=c$. Donc $bd=cd$, ce qui force
$bd=cd=j^2$ (car sinon $j$ serait une valeur propre de $g$ avec multiplicit\'e $\ge 4$).
Toujours d'apr\`es~\eqref{eq:vp}, $(ad,bc)=(j,j^2)$ ou $(ad,bc)=(j^2,j)$.

Si $(ad,bc)=(j,j^2)$, alors $b=c=d$ et $b^2=j^2$, donc $b=c=d=\eta j$ et $a=\eta$, pour un $\eta \in \{\pm 1\}$.
Quitte \`a remplacer $h$ par $-h$, les valeurs propres de $h$ sont donc
$1$, $j$, $j$, $j$ et donc $h$ est bien le seul \'el\'ement
de $\L^{-1}(\{g,g^{-1}\})$ admettant cette liste de valeurs propres.

Si $(ad,bc)=(j^2,j)$, alors $d=jb=jc$ et $b^2=j$, ce qui force $a=b=c=\eta j^2$ et $d=\eta$, pour
un $\eta \in \{\pm 1\}$. Quitte \`a remplacer $h$ par $-h$, les valeurs propres de
$h$ sont donc $1$, $j^2$, $j^2$, $j^2$ et donc $h^{-1}$ est bien le seul \'el\'ement
de $\L^{-1}(\{g,g^{-1}\})$ admettant $1$, $j$, $j$, $j$ comme liste de valeurs propres.
\end{proof}

\bigskip

\begin{coro}\label{coro:conj}
Soit $g \in \Sb\Ob(\wedge^2 V)$ admettant pour valeurs propres $j$, $j$, $j$, $j^2$, $j^2$ et $j^2$.
Alors $g$ et $g^{-1}$ ne sont pas conjugu\'es dans $\Sb\Ob(\wedge^2 V)$.
\end{coro}

\bigskip

\begin{proof}
Supposons trouv\'e $w \in \Sb\Ob(\wedge^2 V)$ tel que $wgw^{-1}=g^{-1}$. Notons $\o$ et $\g$
des ant\'ec\'edents respectifs
de $w$ et $g$ dans $\Sb\Lb(V)$. Alors $\L(\o \g \o^{-1})=g^{-1}=\L(\g^{-1})$.
Cela montre que $\g^{-1}$ est conjugu\'e \`a $\g$ ou \`a $-\g$, ce qui est
impossible suite \`a l'examen des possibles listes de valeurs propres de $\g$
effectu\'e dans la preuve du lemme~\ref{lem:j}.
\end{proof}

\bigskip
\def\jSO{\mathbf{jSO}(\wedge^2 V)}
\def\cubic{\sqrt[3]{\Sb\Lb}(V)}

\section{Construction de ${\boldsymbol{G_{32}}}$}

\medskip

Voyons le groupe $E_6$ comme un sous-groupe fini de $\Ob(\wedge^2 V)$ engendr\'e par des r\'eflexions.
Posons
$$\jSO=\langle \Sb\Ob(\wedge^2 V), j \Id_{\wedge^2 V} \rangle$$
$$\cubic=\{g \in \Gb\Lb(V)~|~\det(g)^3 = 1\} =\langle \Sb\Lb(V),j\Id_V \rangle.\leqno{\text{et}}$$
Alors
\equat\label{eq:jso}
\L^{-1}\bigl(\jSO\bigr)=\cubic.
\endequat
Nous d\'efinissons alors
$$W=\L^{-1}(\langle \Drm(E_6),j\Id_{\wedge^2 V} \rangle).$$
C'est un sous-groupe de $\cubic$. Le but de cette note est de montrer que $W$ est
isomorphe au groupe de r\'eflexions complexes $G_{32}$ de Shephard-Todd. Remarquons
pour commencer que
\equat\label{eq:mu6}
\mub_6 \subset W.
\endequat

\bigskip

\subsection{Reflexions dans ${\boldsymbol{W}}$}
La liste des degr\'es de $E_6$ est $2$, $5$, $6$, $8$, $9$, $12$
alors que celle de ses codegr\'es est $0$, $3$, $4$, $6$, $7$, $10$
(voir~\cite[table~A.3]{broue}).
En particulier, exactement $3$ des degr\'es sont divisibles par $3$, ce qui montre~\cite[th\'eo.~3.4]{springer}
que $E_6$ contient un \'el\'ement $w_3$ admettant la valeur propre $j$ avec multiplicit\'e $3$. Nous noterons
$C_3$ la classe de conjugaison de $w_3$ dans $W$. Puisqu'aussi exactement
$3$ des codegr\'es sont divisibles par $3$, cela implique, par exemple gr\^ace
\`a~\cite[th\'eo.~1.2]{lehrer michel},
que $w_3$ est r\'egulier au sens de Springer~\cite[\S{4}]{springer} (c'est-\`a-dire admet un
vecteur propre pour la valeur propre $j$ dont le stabilisateur dans $E_6$ est trivial). Puisque $E_6$ est
un groupe rationnel, $w_3$ admet aussi $j^2$ comme valeur propre avec multiplicit\'e $3$.
En conclusion, les valeurs propres de $w_3$ sont $j$, $j$, $j$, $j^2$, $j^2$ et $j^2$.
En particulier, $\det w_3=1$ et donc
$$w_3 \in \Drm(E_6) = E_6 \cap \Sb\Ob(\wedge^2 V).$$
D'apr\`es~\cite[th\'eo.~4.2(iii)]{springer},
le centralisateur de $w_3$ dans $E_6$ est
d'ordre $6 \cdot 9 \cdot 12$ ce qui montre que
\equat\label{eq:c3}
|C_3|=2\cdot 5 \cdot 8=80.
\endequat
De plus~\cite[th\'eo.~4.2(iii)]{springer}, si $w \in E_6$, alors
\equat\label{eq:multiplicite}
\text{\it $w \in C_3$ si et seulement si $\dim \Ker(w-j\Id_V)=3$.}
\endequat
Ainsi, si $w \in C_3$, alors $w^{-1} \in C_3$ mais $w^{-1}$ n'est pas conjugu\'e \`a
$w$ dans $\Drm(E_6)$ (en vertu du corollaire~\ref{coro:conj}).

Notons maintenant $\REF(W)$ l'ensemble des r\'eflexions de $W$. Si $s \in \REF(W)$,
les valeurs propres de $\L(s)$ sont $1$, $1$, $1$, $\det(s)$, $\det(s)$, $\det(s)$.
Puisque $\det(s) \in \{j,j^2\}$, les valeurs propres de $\det(s)\L(s)$ sont donc $j$, $j$, $j$,
$j^2$, $j^2$, $j^2$. Il d\'ecoule de~\eqref{eq:multiplicite} que
$\det(s)\L(s) \in C_3$. On a donc d\'efini une application
$$\fonction{\l}{\REF(W)}{C_3}{s}{\det(s)\L(s).}$$
Le r\'esultat suivant nous sera fort utile~:

\bigskip

\begin{lem}\label{lem:ref}
L'application $\l$ est bijective.
\end{lem}

\bigskip

\begin{proof}
Tout d'abord, si $\l(s)=\l(s')$, alors il existe $\xi \in \CM^\times$ tel que
$s'=\xi s$. Puisque $s$ et $s'$ sont des r\'eflexions, cela n'est possible que
si $\xi=1$, et donc $s=s'$. Cela montre l'injectivit\'e de $\l$.

\medskip

Montrons maintenant la surjectivit\'e. Soit $w \in C_3$. D'apr\`es le corollaire~\ref{coro:conj},
il existe un unique $\wti \in \L^{-1}(\{w,w^{-1}\})$ ayant pour valeurs propres
$1$, $j$, $j$, $j$. Alors $j^2\wti$ et $j\wti^{-1}$ sont des r\'eflexions v\'erifiant
$$\l(j^2\wti)=\det(j^2\wti) \L(j^2 \wti)=j^8 \cdot j^4\L(\wti)=\L(\wti)$$
$$\l(j\wti^{-1})=\det(j\wti^{-1})\L(j\wti^{-1})=j^4 \cdot j^2 \L(\wti^{-1})=\L(\wti)^{-1}.\leqno{\text{et}}$$
Donc $w \in \{\l(j^2\wti),\l(j\wti^{-1})\}$, ce qui montre la surjectivit\'e de $\l$.
\end{proof}

\bigskip

On d\'eduit alors de~\eqref{eq:c3} et du lemme~\ref{lem:ref} que
\equat\label{eq:ref}
|\REF(W)|=80.
\endequat

\medskip

\subsection{Structure de ${\boldsymbol{W}}$}
Notre r\'esultat principal est le th\'eor\`eme~\ref{theo:g32} ci-dessous~:
la preuve que nous proposons ici
n'utilise aucune propri\'et\'e du groupe $G_{32}$ ni la classification
des groupes de r\'eflexions complexes et peut donc \^etre vu comme
une construction alternative de $G_{32}$ \`a partir de $E_6$
(nous utilisons en revanche des propri\'et\'es de $E_6$).

\bigskip

\begin{theo}\label{theo:g32}
Le groupe $W$~:
\begin{itemize}
\itemth{a} est d'ordre $155\,520$~;

\itemth{b} est engendr\'e par des r\'eflexions d'ordre $3$~;

\itemth{c} est irr\'eductible et primitif~;

\itemth{d} admet $12$, $18$, $24$, $30$ comme liste de degr\'es.
\end{itemize}
\end{theo}

\bigskip

\begin{proof}
Notons $E_6^\#=\langle \Drm(E_6),j\Id_{\wedge^2 V} \rangle$ et
$W^+=W \cap \Sb\Lb(V)$. Rappelons que $|E_6|=51\,840$. Donc
\equat\label{eq:ordres}
|\Drm(E_6)|=25\,920,\quad |E_6^\#|=77\,760, \quad |W|=155\,520\quad
\text{et}\quad |W^+|=51\,840.
\endequat
Cela montre entre autres~(a). De plus,
\equat\label{eq:centre}
\Zrm(W) =\mub_6 \qquad \text{et}\qquad W/\mub_6 \simeq \Drm(E_6).
\endequat
Soit $\RC=\{\det(s)^{-1}s~|~s \in \REF(W)\}$. Posons
$$G=\langle \REF(W) \rangle\qquad\text{et}\qquad H= \langle \RC \rangle.$$
L'\'enonc\'e~(b) est \'equivalent \`a l'\'enonc\'e
$$W=G.\leqno{(\#)}$$
Tout d'abord, $\L(\RC)=C_3$ et donc $\L(H)=\Drm(E_6)$ (car ce dernier est simple et $C_3$
est une classe de conjugaison). En particulier, $W=H \cdot \mub_6$ et donc,
puisque $H \subset G \cdot \mub_3$, on obtient
$$W=G \cdot \mub_6,$$
ce qui est presque le r\'esultat attendu~$(\#)$. Avant de montrer~$(\#)$, remarquons que,
puisque $\L(H)=\Drm(E_6)$, cela implique que $H$ (et donc $G$) agit de fa\c{c}on
irr\'eductible sur $V$. De plus, si $G$ n'\'est pas primitif, alors $G$ (et donc $H$)
serait monomial~\cite[lem.~2.12]{lehrer taylor}, ce qui impliquerait que
$\L(H)=\Drm(E_6)$ serait aussi monomial, ce qui n'est pas le cas. Donc
$$\text{\it $G$ est irr\'eductible et primitif.}\leqno{(\clubsuit)}$$
L'\'enonc\'e~(c) en d\'ecoule.

\medskip

Nous terminons en montrant simultan\'ement~(b) et~(d).
Soient $d_1$, $d_2$, $d_3$, $d_4$ les degr\'es de $G$.
Il d\'ecoule de~\eqref{eq:ref} et par exemple de~\cite[th\'eo.~4.1]{broue} que
$$
\begin{cases}
d_1d_2d_3d_4=|G|,\\
d_1+d_2+d_3+d_4=84.
\end{cases}\leqno{(\diamondsuit)}
$$
Le morphisme $\det : G \to \mub_3$ est surjectif (car $\det(s) \in \{j,j^2\}$
pour tout $s \in \REF(W)$), ce qui implique que $\mub_3 \subset G$
(car $G/(G \cap \mub_6) \simeq \Drm(E_6)$ est simple). Il reste \`a montrer
que $|G| \neq |W|/2=77\,760$.

Supposons donc que $|G|=77\,760$.
Puisque $\mub_3 \subset G$, tous les $d_i$ sont divisibles par $3$ et, puisque
$\mub_6 \not\subset G$, au moins un des $d_i$ (mettons $d_4$) n'est pas divisible par $6$.
\'Ecrivons $e_i=d_i/3$. Alors
$$
\begin{cases}
e_1e_2e_3e_4=960=2^6\cdot 3\cdot 5,\\
e_1+e_2+e_3+e_4=28,\\
\text{$e_4$ est impair.}
\end{cases}
$$
De la deuxi\`eme \'egalit\'e, on d\'eduit qu'un autre des $e_i$ (mettons $e_3$) est impair.
D'autre part, d'apr\`es la deuxi\`eme \'egalit\'e, au moins l'un des $e_i$ (mettons $e_2$)
est pair et donc $e_1$ est aussi pair. La premi\`ere \'egalit\'e montre que $e_1$ ou $e_2$
(mettons $e_2$) est divisible par $8$. Donc $e_2 \in \{8,16,24\}$.
Une rapide inspection des cas possibles montre qu'alors $e_1=4$, $e_2=16$ et
$\{e_3,e_4\}=\{3,5\}$. Les degr\'es de $G$ sont donc $9$, $12$, $15$ et $48$. Puisque
$16$ divise l'un des degr\'es de $G$, il d\'ecoule alors de~\cite[th\'eo.~3.4(i)]{springer}
que $G=H \times \mub_3$ contient un \'el\'ement d'ordre $16$, donc
$H$ contient un \'el\'ement d'ordre $16$ et donc $\Drm(E_6)=\L(H)$ contient un
\'el\'ement d'ordre $8$, ce qui n'est pas le cas (voir la remarque~\ref{rem:8} ci-dessous pour
une preuve de ce fait n'utilisant que la th\'eorie de Springer).
Cela contredit le fait que $|G|=77\,760$. On a
donc bien montr\'e~$(\#)$, c'est-\`a-dire
$$W=\langle \REF(W) \rangle = G.$$
C'est l'\'enonc\'e~(b).
En particulier, $\mub_6 \subset W$ et donc tous les $d_i$ sont divisibles par $6$.
Posons $a_i=d_i/6$. Alors~($\diamondsuit$) implique que
$$
\begin{cases}
a_1a_2a_3a_4=120=2^3\cdot 3\cdot 5,\\
a_1+a_2+a_3+a_4=14,
\end{cases}
$$
De la m\^eme mani\`ere que pr\'ec\'edemment, puisque $\mub_{12} \not\subset W$,
on peut supposer que $a_3$ et $a_4$ sont impairs et que $a_1$ et $a_2$ sont pairs.
Une rapide inspection des cas possibles implique que $\{a_1,a_2\}=\{2,4\}$ et $\{e_3,e_4\}=\{3,5\}$.
Cela termine la preuve de~(d).
\end{proof}

\bigskip

\begin{coro}\label{coro:g32}
Le groupe $W$ est isomorphe au groupe de r\'eflexions $G_{32}$ de Shephard-Todd.
\end{coro}

\bigskip

\begin{proof}
Cela r\'esulte du th\'eor\`eme~\ref{theo:g32} et de la classification des groupes
de r\'eflexions complexes~\cite{ST}.
\end{proof}

\bigskip

Le corollaire~\ref{coro:g32} donne ainsi une explication
au fait, mentionn\'e dans l'introduction,
que $\Drm(G_{32})/\mub_2 \simeq \Drm(E_6)$~: l'isomorphisme est r\'ealis\'e
par $\L$.

\bigskip

\begin{rema}\label{rem:8}
Dans la preuve du th\'eor\`eme~\ref{theo:g32}, nous avons utilis\'e le fait que $\Drm(E_6)$
ne contient pas d'\'el\'ement d'ordre $8$. Ce fait peut s'obtenir facilement par un calcul
par ordinateur par exemple, mais nous en proposons ici une preuve n'utilisant que la th\'eorie
de Springer. Soit $w \in E_6$ un \'el\'ement d'ordre $8$. Alors $w$ admet forc\'ement parmi
ses valeurs propres une racine primitive huiti\`eme de l'unit\'e $\z$. Comme un seul des
degr\'es de $E_6$ et un seul des codegr\'es de $E_6$ est divisible par $8$, cela implique
que $w$ est un \'el\'ement r\'egulier au sens de Springer~\cite[th\'eo.~1.2]{lehrer michel}.
Alors~\cite[th\'eo.~4.2(v)]{springer}
la liste des valeurs propres de $w$ est $\z^{-1}$, $\z^{-4}$, $\z^{-5}$, $\z^{-7}$, $\z^{-8}$, $\z^{-11}$,
et donc $\det(w)=\z^{-36}=-1$. Donc $w \not\in \Drm(E_6)$.
\end{rema}

\newpage

\def\titre{A construction of the Shephard-Todd group ${\boldsymbol{G_{32}}}$ \\
through the Weyl group of type ${\boldsymbol{E_6}}$}

\def\auteur{{\sc C\'edric Bonnaf\'e}}
%\address{IMAG, Universit\'e de Montpellier, CNRS, Montpellier, France}

\makeatletter
%\email{cedric.bonnafe@umontpellier.fr}
\makeatother

\date{\today}

%\thanks{The author is partly supported by the ANR:
%Projects No ANR-16-CE40-0010-01 (GeRepMod) and ANR-18-CE40-0024-02 (CATORE).}

%

%\tableofcontents

% \noindent{\bf Acknowledgements.} We thank warmly D. Juteau
% for providing us the main ideas leading to Example~\ref{ex:b2-g2}.
%
% \bigskip

\MakeTitle

\newtheorem{theorem}{Theorem}[section]
\def\thetheo{\thesection.\arabic{theorem}}
\newtheorem{proposition}[theorem]{Proposition}
\def\theprop{\thesection.\arabic{proposition}}
\newtheorem{lemma}[theorem]{Lemma}
\def\thelem{\thesection.\arabic{lemma}}
\newtheorem{corollary}[theorem]{Corollary}
\def\thecoro{\thesection.\arabic{corollary}}
\theoremstyle{remark}
\newtheorem{remar}[theorem]{Remark}
\def\thecoro{\thesection.\arabic{remar}}
\renewcommand\theequation{\thetheorem}
\def\equat{\refstepcounter{theorem}\begin{equation}}
\def\endequat{\end{equation}}

\def\abstractname{Abstract}
\def\refname{References}

\begin{abstract}
It is well-known that the quotient of the derived subgroup of the Shephard-Todd
complex reflection group $G_{32}$ (which has rank $4$) by its center is isomorphic to the
derived subgroup of the Weyl group of type $E_6$. We show that this isomorphism can be realized
through the second exterior power, and take the opportunity to propose an alternative construction
of the group $G_{32}$.
\end{abstract}

\pagestyle{myheadings}

\markboth{\sc C. Bonnaf\'e}{\sc A construction of $G_{32}$ through $E_6$}

\bigskip

Let $G_{32}$ denote the complex reflection group constructed by Shephard-Todd~\cite{ST}
and let $E_6$ be a Weyl group of type $\Erm_6$ (which is denoted by $G_{35}$ in the Shephard-Todd
classification). Let $\Sb\pb_4(\FM_{\! 3})$ (resp. $\Sb\Ob_5(\FM_{\! 3})$) denote
the symplectic (resp. orthogonal) group of dimension $4$ (resp. $5$) over the finite field with three
\'el\'ements $\FM_{\! 3}$. Let $\O_5(\FM_{\! 3})$ be the image of $\Sb\pb_4(\FM_{\! 3})$ in
$\Sb\Ob_5(\FM_{\! 3})$ through the natural morphism
$\Sb\pb_4(\FM_{\! 3}) \to \Sb\Ob_5(\FM_{\! 3})$: this is the normal
subgroup of index $2$ of $\Sb\Ob_5(\FM_{\! 3})$.
Finally, if $G$ is a group, let $\Drm(G)$ and $\Zrm(G)$ denote respectively
its derived subgroup and its center and, if $d$ is a non-zero natural number,
let $\mub_d$ be the group of $d$-th roots of unity in $\CM$.

It is shown in~\cite[th\'eo.~8.43~et~8.54]{lehrer taylor}
that $G_{32} \simeq \mub_3 \times \Sb\pb_4(\FM_{\! 3})$
and that $E_6 \simeq \Sb\Ob_5(\FM_{\! 3})$. In particular,
$$\Drm(G_{32})/\mub_2 \simeq \O_5(\FM_{\! 3}) \simeq \Drm(E_6).\leqno{(*)}$$
The purpose of this note is to present a direct elementary explanation
of the isomorphism $\Drm(G_{32})/\mub_2 \simeq \Drm(E_6)$, which in fact allows us to construct the
the complex reflection group $G_{32}$ from the rational reflection group $E_6$.
This construction uses the classical morphism $\Sb\Lb_4(\CM) \longto \Sb\Ob_6(\CM)$, and follows
the same lines as our previous paper~\cite{g31} (in which we constructed
the complex reflection group $G_{31}$ from the Weyl group of type $B_6$).

This note does not pretend to prove a deep result: it is just a
a nice example of the application of the classical theory of reflection groups
(invariants, Springer theory,...).

\bigskip

\noindent{\bf Remark.} It is shown in~\cite[Theo.~8.53]{lehrer taylor} that
$G_{33} \simeq \mub_2 \times \O_5(\FM_{\! 3})$. In particular, this gives an indirect
isomorphism $\Drm(G_{33}) \simeq \Drm(E_6)$. However, we do not know of any construction
of this isomorphism which would be in the same spirit as above.

\bigskip

\section{The morphism $\Sb\Lb_4(\CM) \to \Sb\Ob_6(\CM)$}

\medskip

We recall here the construction of this morphism, using some notation from~\cite{g31}.
Let us fix a complex vector space $V$ of dimension $4$ and let
$$\fonction{\L}{\Gb\Lb(V)}{\Gb\Lb(\wedge^2 V)}{g}{\wedge^2 g}$$
be the natural morphism of algebraic groups. Note that $\wedge^2 V$ has dimension $6$
and that
\equat\label{eq:noyau}
\Ker \L \simeq \mub_2 = \{\pm \Id_V\}.
\endequat
The lists of eigenvalues of elements of $\Gb\Lb(V)$ or $\Gb\Lb(\wedge^2 V)$
will always be given with multiplicities.
If $g \in \Gb\Lb(V)$ admits $a$, $b$, $c$ and $d$ as eigenvalues, then
\equat\label{eq:vp}
\text{\it $\L(g)$ admits $ab$, $ac$, $ad$, $bc$, $bd$ and $cd$ as eigenvalues.}
\endequat
In particular,
\equat\label{eq:det}
\det \L(g) = (\det g)^3.
\endequat

Let us fix now a generator $\e$ of the one-dimensional vector space $\wedge^4 V$.
The choice of this generator allows to identify $\CM$ and $\wedge^4 V$ and to define a
bilinear form
$$\fonction{\b_\wedge}{\wedge^2 V \times \wedge^2 V}{\CM}{(x,y)}{x \wedge y.}$$
This bilinear form is symmetric and non-degenerate. By definition of the determinant,
\equat\label{eq:beta}
\b_\wedge(\L(g)(x),\L(g)(y))=(\det g) \b_\wedge(x,y)
\endequat
for all $g \in \Gb\Lb(V)$ and $x$, $y \in \wedge^2 V$. For dimension and connectedness reasons, the image of
$\Gb\Lb(V)$ through $\L$ is the neutral component $\Cb\Ob(\wedge^2 V)^\circ$
of the conformal orthogonal group
$\Cb\Ob(\wedge^2 V)=\Cb\Ob(\wedge^2 V, \b_\wedge)$ and $\L$ induces an isomorphism of algebraic groups
\equat\label{eq:iso}
\Sb\Lb(V)/\mub_2 \simeq \Sb\Ob(\wedge^2 V).
\endequat

The next elementary lemma will be useful (in this paper, we denote by $j$
a primitive third root of unity).

\bigskip

\begin{lemma}\label{lem:j}
Let $g \in \Sb\Ob(\wedge^2 V)$ having $j$, $j$, $j$, $j^2$, $j^2$ and $j^2$ as eigenvalues.
Then there exists a unique element $\gti \in \Sb\Lb(V)$, having
$1$, $j$, $j$ and $j$ as eigenvalues, and such that
$$\L^{-1}(\{g,g^{-1}\}) = \{\pm \gti,\pm \gti^{-1}\}.$$
\end{lemma}

\bigskip

\begin{proof}
Let $h \in \Sb\Lb(V)$ be such that $\L(h)=g$. Then $\L^{-1}(\{g,g^{-1}\}) = \{\pm h,\pm h^{-1}\}$.
Let $a$, $b$, $c$ and $d$ be the eigenvalues of $h$.
By~\eqref{eq:vp}, by reordering if necessary the eigenvalues of $h$,
we may assume that $ab=ac=j$. In particular, $b=c$. So $bd=cd$, which implies that
$bd=cd=j^2$ (for otherwise $j$ would be an eigenvalue of $g$ with multiplicity $\ge 4$).
Still by~\eqref{eq:vp}, $(ad,bc)=(j,j^2)$ or $(ad,bc)=(j^2,j)$.

If $(ad,bc)=(j,j^2)$, then $b=c=d$ and $b^2=j^2$, so $b=c=d=\eta j$ and $a=\eta$, for some $\eta \in \{\pm 1\}$.
By replacing $h$ by $-h$ if necessary, the eigenvalues of $h$ are then
$1$, $j$, $j$, $j$ and so $h$ is indeed the unique element
of $\L^{-1}(\{g,g^{-1}\})$ admitting this list of eigenvalues.

If $(ad,bc)=(j^2,j)$, then $d=jb=jc$ and $b^2=j$, which implies that $a=b=c=\eta j^2$ and $d=\eta$, for
some $\eta \in \{\pm 1\}$. By replacing $h$ by $-h$ if necessary, the eigenvalues of
$h$ are then $1$, $j^2$, $j^2$, $j^2$ and so $h^{-1}$ is indeed the unique element of
$\L^{-1}(\{g,g^{-1}\})$ admitting $1$, $j$, $j$, $j$ as list of eigenvalues.
\end{proof}

\bigskip

\begin{corollary}\label{coro:conj}
Let $g \in \Sb\Ob(\wedge^2 V)$ having $j$, $j$, $j$, $j^2$, $j^2$ and $j^2$ as eigenvalues.
Then $g$ and $g^{-1}$ are not conjugate in $\Sb\Ob(\wedge^2 V)$.
\end{corollary}

\bigskip

\begin{proof}
Assume that we have found $w \in \Sb\Ob(\wedge^2 V)$ such that $wgw^{-1}=g^{-1}$. Let $\o$ and $\g$
be respective preimages of $w$ and $g$ in $\Sb\Lb(V)$. Then $\L(\o \g \o^{-1})=g^{-1}=\L(\g^{-1})$.
This shows that $\g^{-1}$ is conjugate to $\g$ or to $-\g$, which is impossible
by examining the possible lists of eigenvalues of $\g$ obtained in the proof of Lemma~\ref{lem:j}.
\end{proof}

\bigskip
\def\jSO{\mathbf{jSO}(\wedge^2 V)}
\def\cubic{\sqrt[3]{\Sb\Lb}(V)}

\section{Construction of ${\boldsymbol{G_{32}}}$}

\medskip

Let us see $E_6$ as a finite subgroup of $\Ob(\wedge^2 V)$ generated by reflections.
Set
$$\jSO=\langle \Sb\Ob(\wedge^2 V), j \Id_{\wedge^2 V} \rangle$$
$$\cubic=\{g \in \Gb\Lb(V)~|~\det(g)^3 = 1\} =\langle \Sb\Lb(V),j\Id_V \rangle.\leqno{\text{and}}$$
Then
\equat\label{eq:jso}
\L^{-1}\bigl(\jSO\bigr)=\cubic.
\endequat
We then define
$$W=\L^{-1}(\langle \Drm(E_6),j\Id_{\wedge^2 V} \rangle).$$
It is a subgroup of $\cubic$. The aim of this note is to show that $W$ is isomorphic to
the complex reflection group $G_{32}$ of Shephard-Todd. Note first that
\equat\label{eq:mu6}
\mub_6 \subset W.
\endequat

\bigskip

\subsection{Reflections in ${\boldsymbol{W}}$}
The list of degrees $E_6$ is $2$, $5$, $6$, $8$, $9$, $12$
while its list of codegrees is $0$, $3$, $4$, $6$, $7$, $10$
(see~\cite[Table~A.3]{broue}).
In particular, exactly $3$ of the degrees are divisible by $3$, which shows~\cite[Theo.~3.4]{springer}
that $E_6$ contains an element $w_3$ admitting the eigenvalue $j$ with multiplicity $3$. We will denote by
$C_3$ the conjugacy class of $w_3$ in $W$. Since also exactly
$3$ of the codegrees are divisible by $3$, this implies, for instance by ~\cite[Theo.~1.2]{lehrer michel},
that $w_3$ is regular in the sense of Springer~\cite[\S{4}]{springer} (that is, admits
an eigenvector for the eigenvalue $j$ whose stabilizer in $E_6$ is trivial). Since $E_6$ is
a rational group, $w_3$ also admits $j^2$ as an eigenvalue with multiplicity $3$.
Hence, the eigenvalues of $w_3$ are $j$, $j$, $j$, $j^2$, $j^2$ and $j^2$.
In particular, $\det w_3=1$ and so
$$w_3 \in \Drm(E_6) = E_6 \cap \Sb\Ob(\wedge^2 V).$$
By~\cite[Theo.~4.2(iii)]{springer},
the centralizer of $w_3$ in $E_6$ has order $6 \cdot 9 \cdot 12$, which shows that
\equat\label{eq:c3}
|C_3|=2\cdot 5 \cdot 8=80.
\endequat
Moreover~\cite[Theo.~4.2(iii)]{springer}, if $w \in E_6$, then
\equat\label{eq:multiplicite}
\text{\it $w \in C_3$ if and only if $\dim \Ker(w-j\Id_V)=3$.}
\endequat
Hence, if $w \in C_3$, then $w^{-1} \in C_3$ but $w^{-1}$ is not conjugate to
$w$ in $\Drm(E_6)$ (by Corollary~\ref{coro:conj}).

Now, let $\REF(W)$ denote the set of reflections of $W$. If $s \in \REF(W)$,
the eigenvalues $\L(s)$ are $1$, $1$, $1$, $\det(s)$, $\det(s)$, $\det(s)$.
Since $\det(s) \in \{j,j^2\}$, the eigenvalues of $\det(s)\L(s)$ are then $j$, $j$, $j$,
$j^2$, $j^2$, $j^2$. It follows from~\eqref{eq:multiplicite} that
$\det(s)\L(s) \in C_3$. This defines a map
$$\fonction{\l}{\REF(W)}{C_3}{s}{\det(s)\L(s).}$$
The next result will be very useful:

\bigskip

\begin{lemma}\label{lem:ref}
The map $\l$ is bijective.
\end{lemma}

\bigskip

\begin{proof}
First, if $\l(s)=\l(s')$, then there exists $\xi \in \CM^\times$ such that
$s'=\xi s$. Since $s$ and $s'$ are reflections, this is possible only if
$\xi=1$, and so $s=s'$. This shows that $\l$ is injective.

\medskip

Let us now show the surjectivity. Let $w \in C_3$. By Corollary~\ref{coro:conj},
there exists a unique $\wti \in \L^{-1}(\{w,w^{-1}\})$ admitting
$1$, $j$, $j$, $j$ as eigenvalues. Then $j^2\wti$ and $j\wti^{-1}$ are reflections satisfying
$$\l(j^2\wti)=\det(j^2\wti) \L(j^2 \wti)=j^8 \cdot j^4\L(\wti)=\L(\wti)$$
$$\l(j\wti^{-1})=\det(j\wti^{-1})\L(j\wti^{-1})=j^4 \cdot j^2 \L(\wti^{-1})=\L(\wti)^{-1}.\leqno{\text{and}}$$
So $w \in \{\l(j^2\wti),\l(j\wti^{-1})\}$, which shows that $\l$ is surjective.
\end{proof}

\bigskip

We then deduce from~\eqref{eq:c3} and Lemma~\ref{lem:ref} that
\equat\label{eq:ref}
|\REF(W)|=80.
\endequat

\medskip

\subsection{Structure of ${\boldsymbol{W}}$}
Our main result is Theorem~\ref{theo:g32} below: the proof we propose here
uses neither known properties of the group $G_{32}$ nor the classification
of complex reflection groups and so might be viewed as an alternative construction of
$G_{32}$ starting from $E_6$
(however, note that we use properties of $E_6$).

\bigskip

\begin{theorem}\label{theo:g32}
The group $W$:
\begin{itemize}
\itemth{a} has order $155\,520$;

\itemth{b} is generated by reflections of order $3$;

\itemth{c} is irreducible and primitive;

\itemth{d} admits $12$, $18$, $24$, $30$ as list of degrees.
\end{itemize}
\end{theorem}

\bigskip

\begin{proof}
Set $E_6^\#=\langle \Drm(E_6),j\Id_{\wedge^2 V} \rangle$ and
$W^+=W \cap \Sb\Lb(V)$. Recall that $|E_6|=51\,840$. So
\equat\label{eq:ordres}
|\Drm(E_6)|=25\,920,\quad |E_6^\#|=77\,760, \quad |W|=155\,520\quad
\text{and}\quad |W^+|=51\,840.
\endequat
This shows~(a). Moreover,
\equat\label{eq:centre}
\Zrm(W) =\mub_6 \qquad \text{and}\qquad W/\mub_6 \simeq \Drm(E_6).
\endequat
Let $\RC=\{\det(s)^{-1}s~|~s \in \REF(W)\}$. Set
$$G=\langle \REF(W) \rangle\qquad\text{and}\qquad H= \langle \RC \rangle.$$
The statement~(b) is equivalent to the following one
$$W=G.\leqno{(\#)}$$
First, $\L(\RC)=C_3$ and so $\L(H)=\Drm(E_6)$ (as this last group is simple and $C_3$
is a conjugacy class). In particular, $W=H \cdot \mub_6$ and so,
since $H \subset G \cdot \mub_3$, we obtain
$$W=G \cdot \mub_6,$$
which is almost the expected result~$(\#)$. Before showing~$(\#)$, note that,
since $\L(H)=\Drm(E_6)$, we have that $H$ (and so $G$) acts irreducibly on $V$.
Moreover, if $G$ is not primitive, then $G$ (and so $H$)
would be monomial~\cite[Lem.~2.12]{lehrer taylor}, which would imply that
$\L(H)=\Drm(E_6)$ is monomial, which is false. So
$$\text{\it $G$ is irreducible and primitive.}\leqno{(\clubsuit)}$$
Statement~(c) follows.

\medskip

Let us conclude this proof by showing simultaneously~(b) et~(d).
Let $d_1$, $d_2$, $d_3$, $d_4$ be the degrees of $G$.
It follows from~\eqref{eq:ref} and for instance from~\cite[Theo.~4.1]{broue} that
$$
\begin{cases}
d_1d_2d_3d_4=|G|,\\
d_1+d_2+d_3+d_4=84.
\end{cases}\leqno{(\diamondsuit)}
$$
The morphism $\det : G \to \mub_3$ is surjective (since $\det(s) \in \{j,j^2\}$
for any $s \in \REF(W)$), and this implies that $\mub_3 \subset G$
(because $G/(G \cap \mub_6) \simeq \Drm(E_6)$ is simple). It remains to show that
$|G| \neq |W|/2=77\,760$.

So assume that $|G|=77\,760$.
Since $\mub_3 \subset G$, all the $d_i$'s are divisible by $3$ and, since
$\mub_6 \not\subset G$, at least one of them (say $d_4$) is not divisible by $6$.
Write $e_i=d_i/3$. Then
$$
\begin{cases}
e_1e_2e_3e_4=960=2^6\cdot 3\cdot 5,\\
e_1+e_2+e_3+e_4=28,\\
\text{$e_4$ is odd.}
\end{cases}
$$
From the second equality, we deduce that at least one of the $e_i$'s (say $e_3$) is odd.
From the first equality, we deduce that at least one of the $e_i$'s (say $e_2$) is even
and so $e_1$ is also even. The first equality shows that $e_1$ or $e_2$ (say $e_2$) is divisible by $8$.
So $e_2 \in \{8,16,24\}$. A quick inspection of the possibilities shows that
$e_1=4$, $e_2=16$ and $\{e_3,e_4\}=\{3,5\}$. The degrees of $G$ are then $9$, $12$, $15$ and $48$.
Since $16$ divides one of the degrees,
it then follows from~\cite[th\'eo.~3.4(i)]{springer} that $G=H \times \mub_3$ contains an element
of order $16$ and so $H$ contains an element of order $16$. Therefore, $\Drm(E_6)=\L(H)$ contains an element
of order $8$, which is impossible (see Remark~\ref{rem:8} below for a proof of this fact
based only on Springer theory).
This contradicts the fact that $|G|=77\,760$. So we have shown~$(\#)$, that is,
$$W=\langle \REF(W) \rangle = G.$$
This is statement~(b).
In particular, $\mub_6 \subset W$ and so all the $d_i$'s are divisibles by $6$.
Set $a_i=d_i/6$. Then~($\diamondsuit$) implies that
$$
\begin{cases}
a_1a_2a_3a_4=120=2^3\cdot 3\cdot 5,\\
a_1+a_2+a_3+a_4=14,
\end{cases}
$$
By the same argument as before, since $\mub_{12} \not\subset W$,
we may assume that $a_3$ and $a_4$ are odd and that $a_1$ and $a_2$ are even.
A quick inspection of the possibilities shows that $\{a_1,a_2\}=\{2,4\}$ and $\{e_3,e_4\}=\{3,5\}$.
This concludes the proof of~(d).
\end{proof}

\bigskip

\begin{corollary}\label{coro:g32}
The group $W$ is isomorphic to the reflection group $G_{32}$ of Shephard-Todd.
\end{corollary}

\bigskip

\begin{proof}
This follows from Theorem~\ref{theo:g32} and from the classification of complex reflection groups~\cite{ST}.
\end{proof}

\bigskip

Hence, Corollary~\ref{coro:g32} gives an explanation for the fact, mentioned in the introduction,
that $\Drm(G_{32})/\mub_2 \simeq \Drm(E_6)$: the isomorphism is realized by $\L$.

\bigskip

\begin{remar}\label{rem:8}
In the proof of Theorem~\ref{theo:g32}, we have used the fact that $\Drm(E_6)$
does not contain any element of order $8$. This fact can be easily obtained by a computer
calculation for instance, but we propose here a proof using only Springer theory. Let $w \in E_6$
be an element of order $8$. Then $w$ necessarily admits an eigenvalue which is a primitive
$8$-th root of unity $\z$. As only one of the degrees of $E_6$ and only one of the codegrees of $E_6$
is divisible by $8$, this implies that $w$ is a regular element in the sense of
Springer~\cite[Theo.~1.2]{lehrer michel}.
Then~\cite[Theo.~4.2(v)]{springer}
the list of eigenvalues of $w$ is $\z^{-1}$, $\z^{-4}$, $\z^{-5}$, $\z^{-7}$, $\z^{-8}$, $\z^{-11}$,
and so $\det(w)=\z^{-36}=-1$. So $w \not\in \Drm(E_6)$.
\end{remar}

\end{document}